# Deformation of the myocardium during CPR

Authors: Jafar Moradicheghamahi[1], Gerard Fortuny[1], Josep M. López[1], Joan Herrero[2], Dolors Puigjaner[1]

[1] Departament d'Enginyeria Informàtica i Matemàtiques, Universitat Rovira i Virgili, Spain
[2] Departament d'Enginyeria Química, Universitat Rovira i Virgili, Spain

**Abstract**

Cardiopulmonary resuscitation (CPR) is an emergency procedure performed on patients during cardiac and respiratory arrest. This procedure externally activates the cardiac and respiratory systems via the delivery of chest compression and artificial ventilation. As the main purpose of CPR is to recirculate the blood flow, prediction of the myocardium behavior has great importance. This prediction allows us to have a better understanding of the needed force to recirculate blood without hurting the heart. Finite element method offer the possibility of noninvasive quantification of myocardial deformation. This method is attractive to use for the assessment of myocardial function.

To investigate the behavior of the heart wall, a 3D model of thoracic organs has been prepared using medical images. In this study, to simulate the behavior of different organs, Code-Aster open software is used. For every organ, the material properties are defined. The most important parameters in the study are displacement, normal stress, and Von-Mises stress in the myocardium. Using these parameters, displacement and stress distribution have been predicted.

Effects of the applied force on the chest during CPR and deformation of the myocardium have been predicted by the finite element model. A linear deformation is observable for each organ during force application. Besides, the final location of the heart and ribs and also involved parameters in predicting myocardium deformation are extracted from the model simulations.

This finite element model enables us to have a good vision of the deformation of the myocardium during CPR. Using this method, it is possible to predict the deformation of every part of the heart, especially right and left ventricles.

**Introduction**

An adult heart beats approximately 100,000 times every day, 30 million times every year, and about 2.5 billion times in a lifetime. [1]. Iin spite of having a small size, it circulates 7 m³ of blood daily, 2500 m³ every year, and 200000 in an average lifetime[2]. Annually, about 275,000 Europeans experience cardiac arrest and receive cardiopulmonary resuscitation (CPR), and only 29,000 discharged from the hospital [3]. For this reason, cardiac arrest is mentioned as the main cause of death. [4]. Cardiopulmonary resuscitation, an immediate procedure for patients with cardiac arrest, has a high potential for recirculation of heart and air in the lungs by compression of the chest and artificial ventilation of the lungs [5]. In recent decades, CPR has become a significant emergency procedure to save the lives of patients with cardiac arrest [6]. Blood flow during CPR is correlated with the compression depth [7]. The main purpose of CPR is to supply oxygenated blood to the myocardium and other organs, to fully restore the function of heart and brain. However, the basic physiological principles for the CPR life-saving process are just partially understood and often argumentative. [4].

The cost of numerical simulation is not high, it can test the system in different conditions, and it is able to predict parameters that are difficult to measure in the laboratory [8]. Numerical modeling of the myocardium is an acceptable tool in understanding the mechanisms of healthy and abnormal behavior of the heart and considering the mechanical activity of tissue is an important feature of this modeling [9-11].

The myocardium is the thick muscular layer of the heart. it comprises muscle fibers, extracellular matrix, blood vessels, blood, and interstitial fluid [12]. In most models used for the left ventricular wall (LV) or other parts of the heart, the epicardium (outer layer) and endocardium (inner layer), which are much thinner than the myocardium, are neglected [13].

Acceptable computational models for left ventricular behavior can be developed using some simplifications in myocardial deformation. In such models, the myocardium is almost represented by a shape that can be defined by several parameters. In most cases, the shapes are symmetrical about the left ventricular axis. Examples of this class of models include cylindrical models [14], spherical models [15], and prolate spheroidal models [16]. With advances in technology, models based on the anatomy of the heart, the mechanism, and interactions of the cardiac system with geometric details are widely used and predict more realistic results. [17-19]. Usually, in these models, finite element method is applied, which by using a great number of freedom degrees, predict heart deformation and other desired variables. [20].

To enhance CPR safety, and have a reduction in occurrences of the relevant injuries, a wider vision of the consequences of physiological parameters related to performing CPR must be provided. Therefore, precise and efficient tools for simulating and evaluating CPR results are needed to make the possibility to develop clinical guidelines and technological inventions. For this, there is a need for models which enable us to have a better understanding of the injury risks [21]. Despite the great wide range of studies on normal heart function and also deformation of the chest during CPR, do date there is no study concentrating on deformation of the myocardium during CPR.

A comprehensive FEM model for chest and other organs provides an appropriate opportunity for the prediction of heart deformation during CPR. For this purpose, in the present study, a 3D dictionary type model from a 3D whole-body model was extracted to investigate the deformation of the heart and stress distribution during CPR. This model includes ribs, intercostal muscles, costal cartilages, sternum and heart.

**Problem statement**

The 3D model used in the simulation has been shown in Figure 1 (a) which includes different organs such as ribs and heart. In Figure 1 (b) detailed view of the heart has been presented. This model of the heart includes right and left atriums, right and left ventricles, and myocardium. In this study, all of the surrounding muscles and tissues were simulated as a "soft" part. Also, the heart was modeled in two parts; myocardium and blood. However, the effect of valves was neglected.

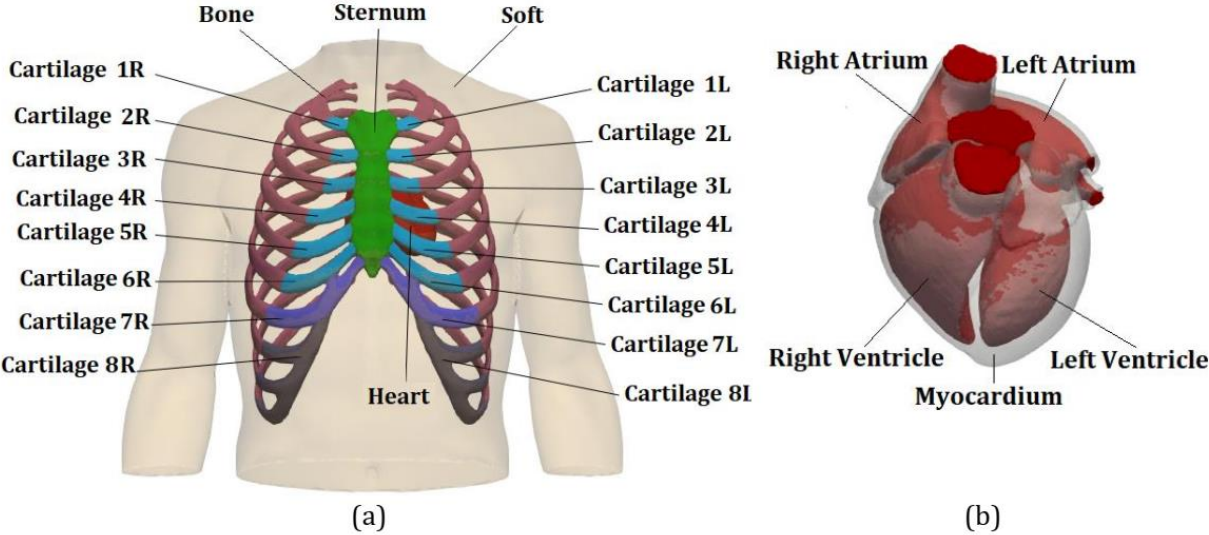

Figure 1: Details of the simulated geometry for CPR

The finite element method was employed to provide a numerical solution for the governing equations. As illustrated in Figure 2, a grid consisting of 499943 nodes and 2937270 elements was used for the analysis. 75260 elements are 2D triangular and 2862010 elements are 3D tetrahedral.

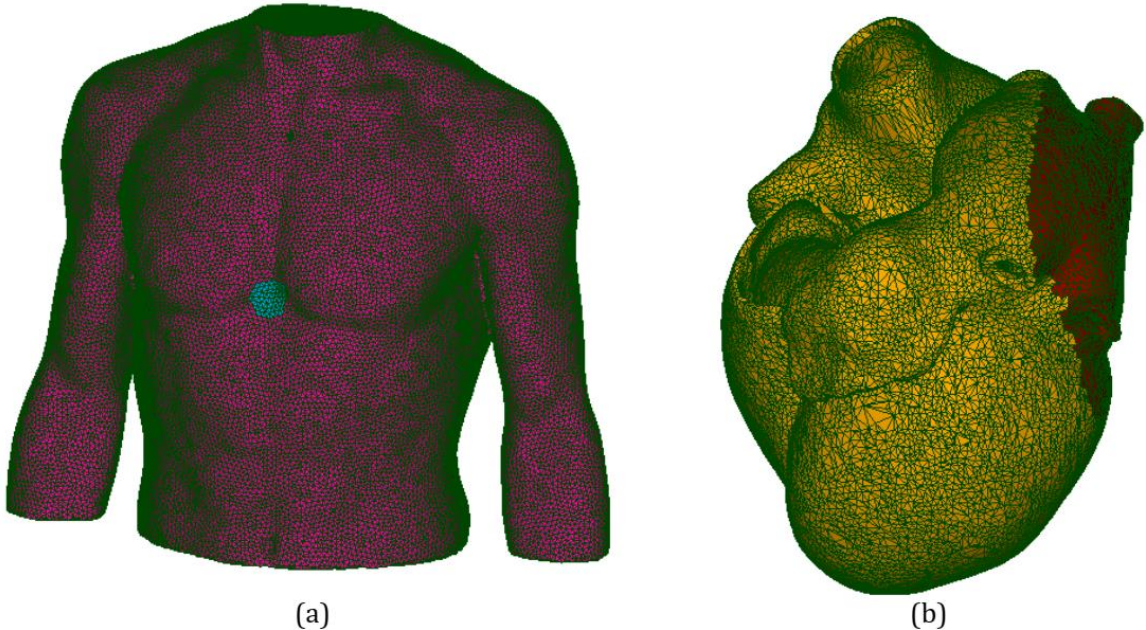

Figure 2: Schematic of the prepared grid for the whole model

In this study, linear elastic models were applied for all materials. Due to its accuracy in predicting small deformations [22], the linear model was preferred to non-linear models which requires a higher computational cost.

As material properties depend on many factors such as age and gender, they are different in various publications. Based on what has been said in [23] Young module was 1e6 Pa for the heart and Poisson's ratio was 0.3 in the present study. Other organs' properties have been shown in Table 1. It should be noted that the Young module value is not the same for different cartilages. The location of cartilages have been shown in Figure 1.

Table 1: Material properties of different parts

| Organ | Young module (Pa) | Poisson ratio | Density (kg/m$^3$) | Reference |
| --- | --- | --- | --- | --- |
| Myocardium | 1e6 | 0.3 | 2000 | [23] |
| Bones | 2e9 | 0.2 | 1000 | [24] |
| Soft tissue | 9e4 | 0.2 | 1000 | [25] |
| Cartilage 1L, 1R | 9e6 | 0.3 | 1100 | [26] |
| Cartilage 2L, 2R, 3L, 3R | 8e7 | 0.3 | 1100 | [26] |
| Cartilage 4L, 4R, 5L, 5R | 7e7 | 0.3 | 1100 | [26] |
| Cartilage, 6L,6R, 7L, 7R, 8L, 8R | 4e7 | 0.3 | 1100 | [26] |

In the present model, the applied force has been defined on the sternum and in the perpendicular direction to the chest. For this purpose, as is illustrated in Figure 3 (a) with blue color, an area was modeled to apply the force. The value of the force was 1e6 N on the applied area of 498.68 mm$^2$. To simulate CPR, some fixed areas are required. In order to this, an areas in the back of the chest were fixed. Also in this model, we fixed a part of the myocardium (blue area in Figure 3 (b)).

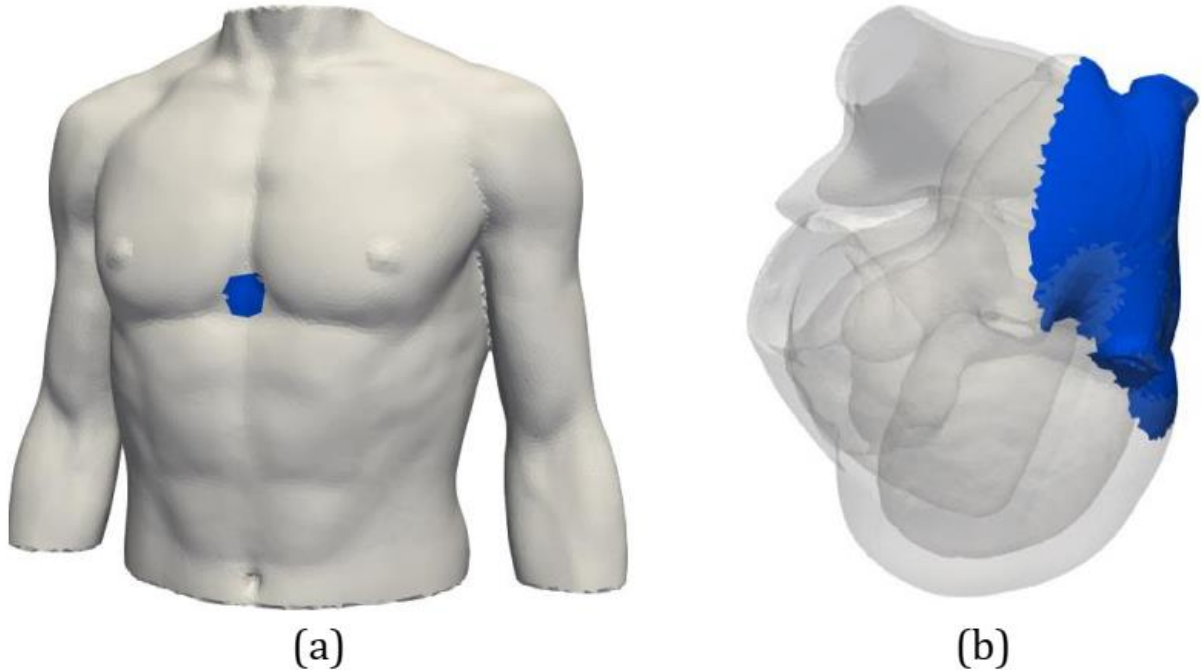

Figure 3; (a): the area of the applied force (b): The fixed area of the myocardium

The simulation was performed for 0.5 s and the defined time step was 0.05 s. In The following, results for four times during CPR have been discussed.

**Results**

In this study, using a model which includes ribs, intercostal muscles, costal cartilages, sternum, heart, and surrounding muscles, the CPR procedure was simulated numerically. As there are several organs in the model, and as the main goal of CPR is recirculating blood in the heart, we focused on the deformation of the myocardium in this study. To this purpose, three parameters were investigated; Displacement, normal stress, and von Mises stress. Results have been presented for four different times including 0.4, 0.45, and 0.5 s.

With the intention of more accurate investigation, in addition to contours shown on the myocardium, some contours related to a cross-section of the myocardium (as shown in Figure 4) have been prepared. This cross-section was chosen to show right and left ventricles.

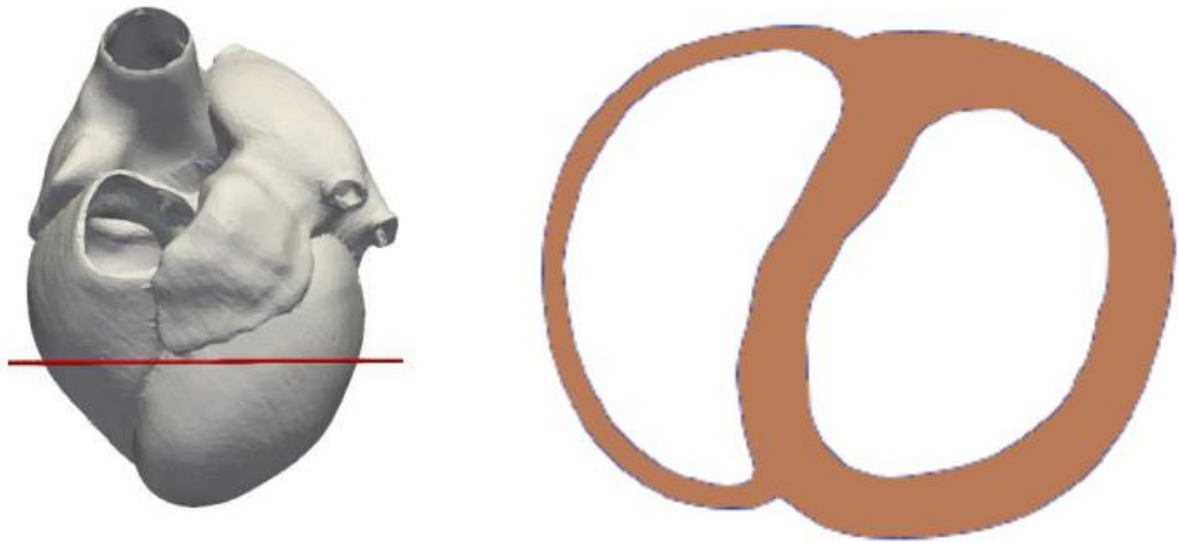

Figure 4: Details of the cross-section

As the right ventricle is the closest part of the heart to the sternum (where the force is applied), has the most deformation compared to other parts of the heart (2.3e -4 m as shown in Figure 5 and Figure 6). In addition to the deformation of the right ventricle, the right atrium deformed considerably. The right parts of the heart must transmit blood to the lungs through pulmonary arteries to oxygenate it. Due to this fact, deformation of the right atrium and right ventricle has high importance. Also, based on the observation in Figure 6 there is a slight deformation in the interventricular septum which can help to flow the blood from the left ventricle to the aorta.

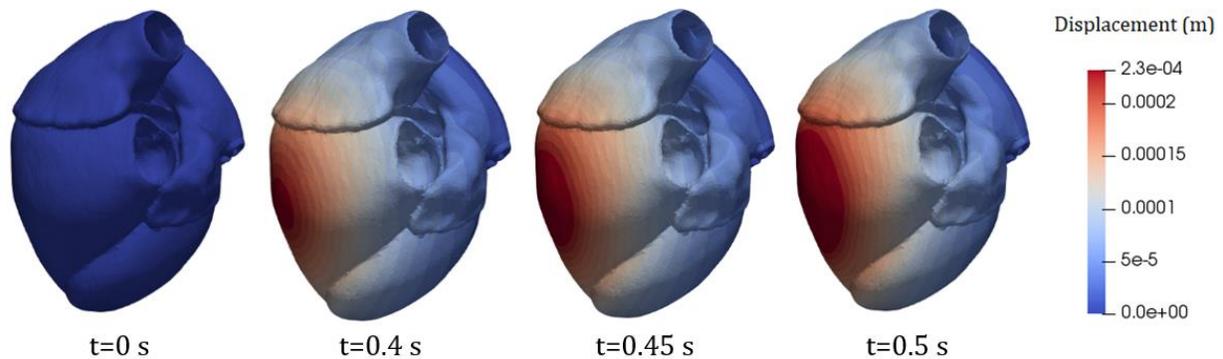

Figure 5: Displacement contours for different times

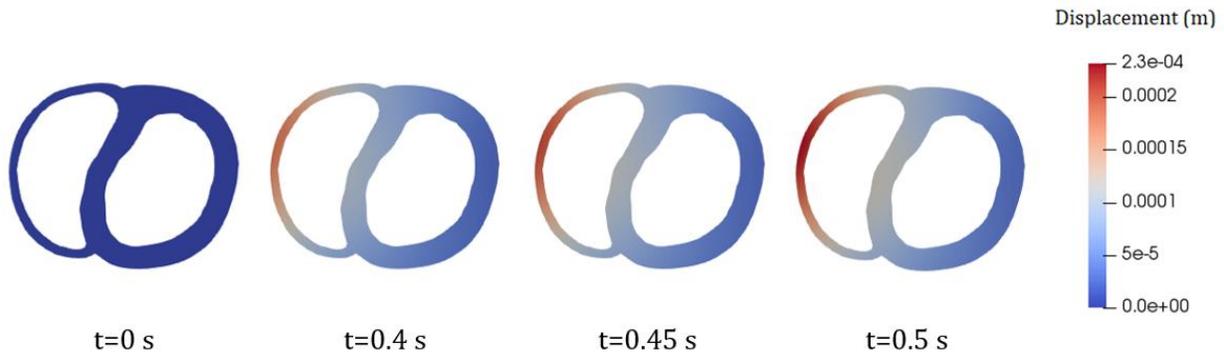

Figure 6: Displacement contours in cross-section for different times

As evident in Figure 5 and Figure 6, due to the fixed area in the left side of the heart, the deformation of the left atrium and left ventricle is not considerable. After entering the right atrium, because of the one direction valves function, blood can move only in one direction. Based on this fact, the deformation of the right atrium and right ventricle is sufficient to recirculate the blood and there is no necessity for the deformation of the left side of the heart.

Normal stress and von Mises stress contours have been presented in Figure 7, Figure 8, and Figure 9. As it is obvious in Figure 7 and at t=0.5 s, the maximum value of the normal stress is in the right ventricle wall (1.3 e +5 Pa). It worth to be noted that values of normal stress are noticeable in the right atrium and left ventricle. By more investigation in Figure 8, it could be concluded that the maximum value of the stress is in the right ventricle and close to the interventricular septum. Another finding based on Figure 8 is the difference between stress values for the inner and outer layers of the myocardium. It is observable that the outer layer of the heart muscle is affected by applied force during CPR significantly more than the inner layer.

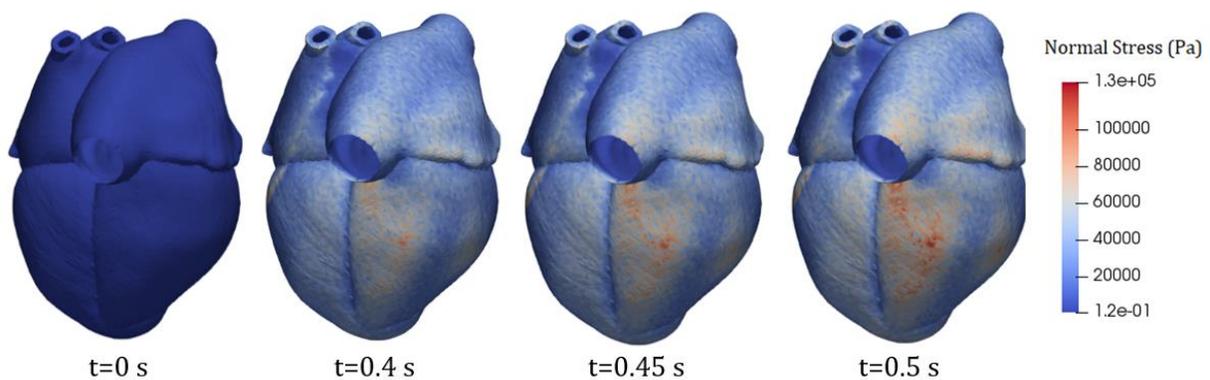

Figure 7: Normal stress distribution in the myocardium for different times

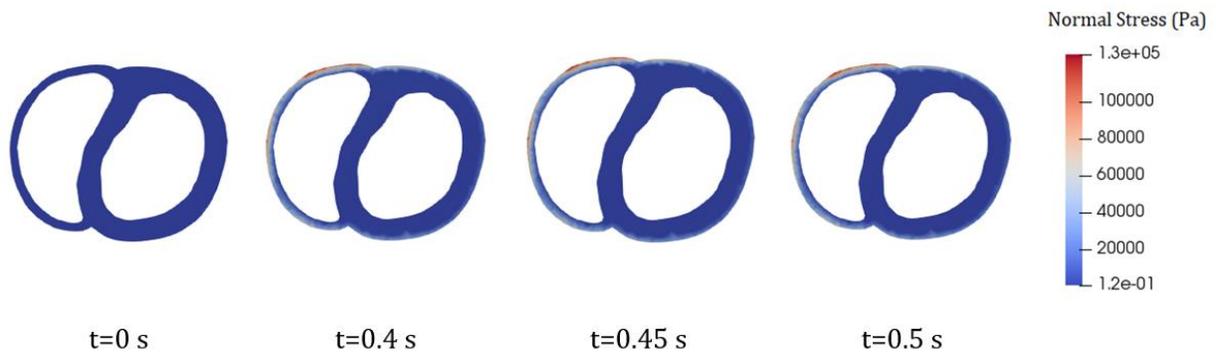

Figure 8: Normal stress distribution in the cross-section for different times

In Figure 9 von Mises stress contours have been presented for different times. Similar to normal stress, the maximum value for von Mises stress is for t=0.5 s and in the right ventricle. In this figure, the influence of the applied force on the other parts of the heart is evident. This influence is such that in the left atrium there are some areas with high values of von Mises stress.

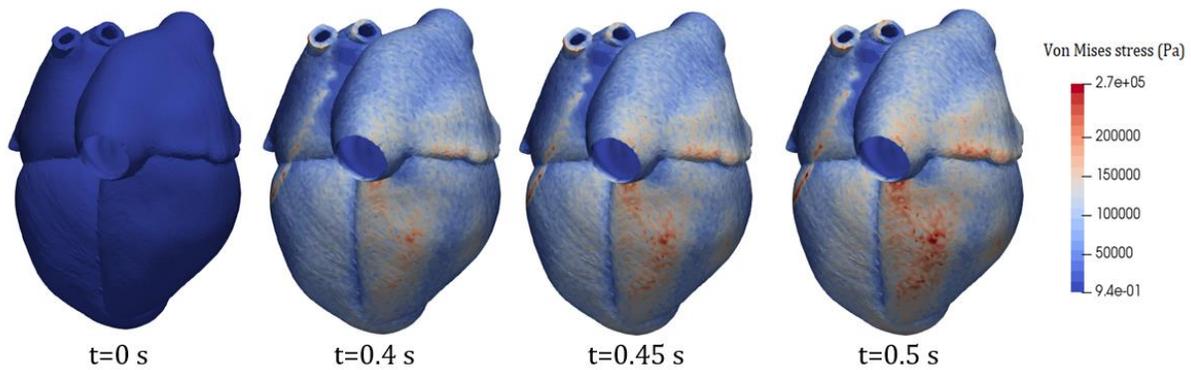

Figure 9: Von Mises stress distribution in the myocardium for different times

**Conclusion**

The main purpose of CPR is to recirculate the blood in the heart. In this study using FEM, the CPR procedure was simulated on a realistic model consisting of different organs in the human chest. It was observed that the maximum displacement is in the right ventricle wall which helps to flow the blood to lungs. Also, displacement of the interventricular septum is notable which is helpful to send blood from the left article to different parts of the body through the aorta. The maximum value of the stress is in the area of the right ventricle near to interventricular septum.

In a conclusion, it can be claimed that this model has a great potential to predict deformation and stress distribution in the heart and other parts of the human chest during CPR.